\documentclass[11pt,onecolumn]{article}
\usepackage{color,graphicx,xcolor,cite}
\usepackage{amsfonts,amssymb,amsmath,amsthm}
\usepackage{empheq}
\usepackage[ruled]{algorithm2e}
\usepackage{url}
\usepackage{multirow,array,booktabs}
\usepackage{threeparttable}
\usepackage{cases} 

\usepackage{CJK}
\usepackage{latexsym}
\usepackage{color}
\definecolor{myblue}{rgb}{0,0,0.5}
\definecolor{mygreen}{rgb}{0,0.5,0}
\definecolor{myred}{rgb}{0.5,0,0}

\usepackage{amsmath,amsthm,amssymb}
\usepackage{graphicx}
\usepackage{hyperref}
\newcommand{\nn}{\nonumber}

\def \[{\begin{equation}}
\def \]{\end{equation}}

\newtheorem{theorem}{Theorem}[section]

\newtheorem{lemma}{Lemma}[section]

\newtheorem{remark}{Remark}[section]

\newtheorem{proposition}{Proposition}[section]

\newcommand{\hdp}[1]{\begin{center}\fbox{
    \begin{minipage}{15.5cm}#1 \end{minipage}}\end{center}\vspace{1ex minus 1ex}}

\begin{document}
\begin{CJK*}{GBK}{song}

\begin{center}

{\large \bf  Balanced Augmented Lagrangian Method for Convex Programming
}\\

\bigskip
\medskip

  {\bf Bingsheng He}\footnote{\parbox[t]{15.5cm}{
 Department of Mathematics,  Nanjing University, Nanjing, China.
  This author was supported by the NSFC Grant 11871029. Email: hebma@nju.edu.cn}}
 \quad
   \quad
 {\bf Xiaoming Yuan}\footnote{\parbox[t]{15.5cm}{Department of Mathematics, The University of Hong Kong, Hong Kong, China. Email:  xmyuan@hku.hk
  }}

\medskip

\today

\end{center}

\bigskip

{\small

\parbox{0.95\hsize}{

\hrule

\medskip

{\bf Abstract.} We consider the convex minimization model with both linear equality and inequality constraints, and reshape the classic augmented Lagrangian method (ALM) by balancing its subproblems. As a result, one of its subproblems decouples the objective function and the coefficient matrix without any extra condition, and the other subproblem becomes a positive definite system of linear equations or a positive definite linear complementary problem. The balanced ALM advances the classic ALM by enlarging its applicable range, balancing its subproblems, and improving its implementation. We also extend our discussion to two-block and multiple-block separable convex programming models, and accordingly design various splitting versions of the balanced ALM for these separable models. Convergence analysis for the balanced ALM and its splitting versions is conducted in the context of variational inequalities through the lens of the classic proximal point algorithm.

\medskip


 \medskip

  \hrule

  }}

\bigskip

\section{Introduction}

The classic augmented Lagrangian method (ALM) was proposed in \cite{Hes69,Powell69}, and since then it has been playing fundamental roles in algorithmic design for various convex programming problems. For instance, it is the root of the alternating direction method of multipliers (ADMM) proposed in \cite{GM75}, which is nowadays a benchmark algorithm used widely in many areas. We refer to, e.g. \cite{Ber82,CP-Acta,FG83,Glow89,Rock76A}, for insightful discussions on the ALM and its wide applications in different areas such as PDEs, optimization, optimal control, image processing, and scientific computing. In particular, it was shown in \cite{Rock76A} that the ALM is an application of the classic proximal point algorithm (PPA) which was originally proposed in \cite{Mar70}.

Let us start with the following canonical convex minimization model with linear equality constraints:
\[  \label{Problem-E}
    \min \{{\theta}(x) \; | \;  Ax =b, \, x\in {\cal X} \},
   \]
where ${\theta}: {\Re}^{n}\to {\Re}$ is a closed proper convex but not necessarily smooth function; ${\cal X} \subseteq \Re^{n}$ is a closed convex set; $A \in \Re^{m\times n}$; and $b\in \Re^{m}$. The iterative scheme of ALM for (\ref{Problem-E}) reads as
\begin{subequations} \label{ALM-E}
\begin{numcases}{\hskip-1.5cm\hbox{(ALM)\quad}}
\label{ALM-E-x} \begin{array}{rcl} x^{k+1} \in \arg\min \bigl\{\theta(x)  -(\lambda^k)^T(Ax-b)+  {\displaystyle \frac{r}{2}}
        \|Ax-b\|^2    \; \big| \; x\in{\cal X} \bigr\} ,
                                                                           \end{array} \\[0.1cm]
\label{ALM-E-y}  \begin{array}{rcl} \lambda^{k+1}= \lambda^k - r (Ax^{k+1} -b),
                                                                           \end{array}
\end{numcases}
\end{subequations}
in which $r>0$ is the penalty parameter and $\lambda \in \Re^m$ is the Lagrange multiplier. Hereafter, $x$ and $\lambda$ are referred to the primal and dual variables, respectively. In general, the subproblem (\ref{ALM-E-x}) needs to be solved iteratively and thus outer-inner nested iterations are rendered to implement the ALM (\ref{ALM-E}). Therefore, how to solve the $x$-subproblem (\ref{ALM-E-x}) determines the difficulty of implementing the ALM (\ref{ALM-E}). An obvious obstacle is that the objective function $\theta(x)$, the coefficient matrix $A$, and the set ${\cal X}$ are all aggregated to be considered simultaneously in the $x$-subproblem (\ref{ALM-E-x}). Thus, the $x$-subproblem (\ref{ALM-E-x}) dominates the computation while the $\lambda$-subproblem (\ref{ALM-E-y}) is trivial. In this sense, these two subproblems in the classic ALM (\ref{ALM-E}) are unbalanced.

In this paper, we suggest to decouple the objective function $\theta(x)$ and the coefficient matrix $A$ in the subproblem (\ref{ALM-E-x}) so as to alleviate this subproblem substantially, and then shift the consideration of the matrix $A$ to the subproblem (\ref{ALM-E-y}). The classic ALM (\ref{ALM-E}) is thus reshaped, and the resulting subproblems are balanced in the sense that the difficulty of the $x$-subproblem only depends on $\theta(x)$ and ${\cal X}$, and that of the $\lambda$-subproblem becomes to depend on $A$. This balancing idea has an immediate advantage when the function $\theta(x)$ has the favorable property that its proximity operator can be represented by a closed-form. That is, the proximity operator of the objective function $\theta(x)$, which is defined by
\begin{equation}\label{proximity}
\hbox{Prox}_\theta^r(x):=\arg\min\bigl\{ \theta(y)+{\displaystyle \frac{r}{2}}
        \|y-x\|^2    \; \big| \; y\in{\Re^n}\bigr\}, \;\forall x \in \Re^n, \;\forall r>0,
\end{equation}
has a closed-form representation. This scenario arises in many applications, especially in contemporary data science domains. We refer to, e.g., \cite{CDS,CR09,RFP}, for some applications whose corresponding function $\theta(x)$ usually prompts sparsity- or low-rank properties of a desired solution and hence can be specified as the $l_1$-norm function (or the nuclear-norm function for the case with matrix variables). Our idea of decoupling $\theta(x)$ and $A$ can be further explained by the following motivation. Ignoring some constant terms, we know that the subproblem (\ref{ALM-E-x}) can be rewritten as
$$
x^{k+1} \in \arg\min \bigl\{\theta(x)  +  {\displaystyle \frac{r}{2}}
        \|Ax-b -\frac{1}{r} \lambda^k\|^2  \; \big| \; x\in{\cal X}\bigr\}.
        $$
For the classic ALM (\ref{ALM-E}), even when $\hbox{Prox}_\theta^r(x)$ can be represented by a closed-form and ${\cal X}=\Re^n$, in general the subproblem (\ref{ALM-E-x}) may still be difficult when the matrix $A$ is not identity. If $\theta(x)$ and $A$ are decoupled and the primeval $x$-subproblem (\ref{ALM-E-x}) is replaced by an easier one in form of
\[   \label{CP-sub}
x^{k+1}=\arg\min \bigl\{\theta(x)  +  {\displaystyle \frac{r}{2}  \|x-q^k\|^2  \; \big| \; x\in{\cal X}}  \bigr\} ,  \]
in which $q^k \in \Re^n$ is a certain constant vector, then the solution of (\ref{CP-sub}) can also be given by the closed-form representation of $\hbox{Prox}_\theta^r(x)$ when ${\cal X}=\Re^n$.

Some existing algorithms in the literature can be applied to (\ref{Problem-E}), and $\theta(x)$ and $A$ can be decoupled in their implementations. For example, as analyzed in \cite{IMA-HMY}, we can consider regularizing the objective function of (\ref{ALM-E-x}) with a proximal term. The resulting proximal version of the ALM (PALM for short) can be written as
\begin{subequations} \label{PALM-E}
\begin{numcases}{\hskip-1.5cm\hbox{(PALM)}}
\label{PALM-E-x} \begin{array}{rcl} x^{k+1} \in \arg\min \bigl\{\theta(x)  -(\lambda^k)^T(Ax-b)+  {\displaystyle \frac{r}{2}}
        \|Ax-b\|^2    + {\displaystyle\frac{1}{2}}\|x-x^k\|_G^2  \; \big| \; x\in{\cal X}  \bigr\},
                                                                           \end{array} \\[0.1cm]
\label{PALM-E-y}  \begin{array}{rcl} \lambda^{k+1}= \lambda^k - r (Ax^{k+1} -b),
                                                                           \end{array}
\end{numcases}
\end{subequations}
with $G \in \Re^{n\times n}$ and the notation $\|x\|_G^2:=x^TGx$. If we choose $G=\sigma I_n-r A^TA$ in (\ref{PALM-E-x}) with $\sigma>0$, then the generic PALM (\ref{PALM-E}) is specified as

\begin{subequations} \label{L1PALM-E}
\begin{numcases}{\hskip-1.5cm\hbox{(LALM)}}
\label{L1PALM-E-x} \begin{array}{rcl} x^{k+1} \in \arg\min \bigl\{\theta(x)  +  {\displaystyle\frac{\sigma}{2}}
        \big\|x-(x^k+{\displaystyle\frac{1}{\sigma}}A^T(\lambda^k-r(Ax^k-b))) \big\|^2 \; \big|\;x\in {\cal X}\bigr\},
                                                                           \end{array} \\[0.1cm]
\label{L1PALM-E-y}  \begin{array}{rcl} \lambda^{k+1}= \lambda^k - r (Ax^{k+1} -b),
                                                                         \end{array}
\end{numcases}
\end{subequations}
in which the subproblem (\ref{L1PALM-E-x}) is in form of (\ref{CP-sub}) and thus it is reduced to $\hbox{Prox}_\theta^r(x)$ when ${\cal X}=\Re^n$. The scheme (\ref{L1PALM-E}) is called the linearized ALM (LALM for short) because the quadratic term $\frac{r}{2} \|Ax-b\|^2$ in (\ref{PALM-E-x}) is ``linearized" by $G=\sigma I_n-r A^TA$. From an analytical point of view, it is easy to see that if $\sigma$ is large enough such that $\sigma >r\|A^TA\|$, then the matrix $G=\sigma I_n-r A^TA$ is positive definite, and essentially convergence of the LALM (\ref{L1PALM-E}) can be conducted by following existing works such as \cite{Esser,He02MP,WY,YY13}. This means the classic ALM (\ref{ALM-E}) can be revised as the LALM (\ref{L1PALM-E}) to decouple $\theta(x)$ and $A$. But the subproblem  (\ref{L1PALM-E-x}) is correlated implicitly with $A$ because of the condition $\sigma >r\|A^TA\|$. On the other hand, note that the approximation of (\ref{L1PALM-E-x}) to the primeval $x$-subproblem (\ref{ALM-E-x}) is less accurate when $\sigma$ is larger, because of the higher weight of the additional quadratic term $\frac{1}{2}\|x-x^k\|_G^2$ with $G=\sigma I_n-r A^TA$. When $\|A^TA\|$ is large, $\sigma$ is forced to be large, and the consequence is that the step size for solving (\ref{L1PALM-E-x}) becomes small and it is doomed that more outer iterations are needed, despite that the inner iterations can be avoided. In \cite{IMA-HMY}, it is shown that the best bound of $\sigma$ is $0.75\cdot r\|A^TA\|$ to ensure the convergence of (\ref{L1PALM-E}), while the mentioned difficult remains if $\|A^TA\|$ is too large.

There is another algorithm that can be applied to the problem (\ref{Problem-E}), and $\theta(x)$ and $A$ can be decoupled in its implementation. More specifically, let us consider the Lagrangian function of (\ref{Problem-E}) and its saddle-point reformulation, and then apply the primal-dual method proposed in \cite{CHPock}. With some tedious details skipped, the resulting iterative scheme can be written as
\begin{subequations} \label{CP-PPA}
\begin{numcases}{}
\label{CP-x} \begin{array}{rcl} x^{k+1} = \arg\min \bigl\{\theta(x)  +  {\displaystyle\frac{r}{2}}
                                                                                 \|x - (x^k +  {\displaystyle\frac{1}{r}} A^T\lambda^k)\|^2    \; \big| \; x\in{\cal X} \bigr\}.
                                                                           \end{array} \\
\label{CP-y}\begin{array}{rcl} \lambda^{k+1}= \lambda^k -  {\displaystyle\frac{1}{s}} \bigl(A(2x^{k+1} - x^k) -b \bigr),
                                                                           \end{array}
\end{numcases}
\end{subequations}
where $r>0$ and $s>0$ are parameters for the primal- and dual-variable subproblems, respectively. In (\ref{CP-x}), $\theta(x)$ and $A$ are also decoupled, and this subproblem is also reduced to $\hbox{Prox}_\theta^r(x)$ when ${\cal X}=\Re^n$. Nearly at the same time as \cite{CHPock}, the primal-dual method proposed in \cite{CHPock} was explained as an application of the classic PPA in \cite{HY-SIIMS}, and then this PPA explanation has been used to analyze the convergence for variants of the primal-dual method (\ref{CP-PPA}), as well as other first-order algorithms, in the literature, see, e.g. \cite{Becker,CHPockMP,PC2001}. To ensure the convergence of (\ref{CP-PPA}), as analyzed in \cite{CHPock}, the condition
\[   \label{Prameter-rs}     rs  > \|A^TA\|  \]
is required. Following the PPA explanation in \cite{HY-SIIMS}, the condition (\ref{Prameter-rs}) is used to ensure the positive definiteness of the matrix that is used to definite the underlying PPA. We refer to, e.g. \cite{CHPock,CP-Acta,HY-SIIMS, COAP-HYZ} for some efficient applications of the primal-dual method (\ref{CP-PPA}) to some image reconstruction problems whose corresponding $\|A^TA\|$ is small. Therefore, despite that the objective function $\theta(x)$ and the coefficient matrix $A$ are decoupled in notation, the subproblem (\ref{CP-x}) is correlated implicitly with $A$ because of the condition (\ref{Prameter-rs}). Clearly, the same difficulties as those for implementing the PALM (\ref{PALM-E}) should be tackled if $\|A^TA\|$ is too large.

Our main purpose is to balance the subproblems of the classic ALM (\ref{ALM-E}) such that both subproblems could be easy for some applications. More specifically, let $r>0$ and $\delta>0$ be arbitrary constants; and define the positive definite matrix $H_0\in \Re^{m \times n}$ as
\[\label{H0}
 H_0:=\bigl(\frac{1}{r}AA^T + \delta I_m  \bigr).
 \]
Then, with $q_0^k:=x^k + \frac{1}{r} A^T\lambda^k$, the classic ALM (\ref{ALM-E}) for the problem (\ref{Problem-E}) is balanced as
\begin{subequations} \label{SD-PPA}
\begin{numcases}{\hskip-1.0cm\hbox{(Balanced ALM)}}
\label{SD-x} \begin{array}{rcl} x^{k+1}    & = & \arg\min \bigl\{\theta(x)  +    {\displaystyle \frac{r}{2}}
                                                                                 \|x - q_0^k\|^2   \; \big| \; x\in{\cal X} \bigr\},
                                                                           \end{array} \\
\label{SD-l}\begin{array}{rcl}
            \lambda^{k+1}        &= & \lambda^k-H_0^{-1}(A(2 x^{k+1}-x^k) - b).                                                                         \end{array}
\end{numcases}
\end{subequations}
Hence, for the model (\ref{Problem-E}) with ${\cal X}=\Re^n$, the balanced ALM (\ref{SD-PPA}) is reduced to
\begin{subequations} \label{SD-PPA-short}
\begin{numcases}{}
\label{SD-x-short} \begin{array}{rcl} x^{k+1}    & = & \hbox{Prox}_\theta^r(q_0^k), \end{array} \\
\label{SD-l-short}\begin{array}{rcl} \lambda^{k+1}        &= & \lambda^k-H_0^{-1}(A(2 x^{k+1}-x^k) - b).     \end{array}
\end{numcases}
\end{subequations}
In the $x$-subproblem (\ref{SD-x}), it is easy to discern that $\theta(x)$ and $A$ are decoupled while the parameter $r$ is not restricted by any condition related to $\|A^TA\|$ explicitly or implicitly, and thus it could be as easy as estimating $\hbox{Prox}_\theta^r$ when ${\cal X}=\Re^n$. Moreover, the $\lambda$-subproblem (\ref{SD-l}) is still very easy though it involves the matrix $A$ and thus becomes slightly more difficult than (\ref{ALM-E-y}) or (\ref{CP-y}), see Remark \ref{easyH} for more details. In this sense, the subproblems of the classic ALM (\ref{ALM-E}) are balanced in (\ref{SD-PPA}). Obviously, the balanced ALM (\ref{SD-PPA}) enjoys the proximity-induced feature while it can avoid possible tiny step sizes for the subproblems (\ref{SD-x}) even when $\|A^TA\|$ is large. This is an essential difference of the balanced ALM (\ref{SD-PPA}) from the PALM (\ref{PALM-E}) and the primal-dual method (\ref{CP-PPA}). We consider the balanced ALM (\ref{SD-PPA}) a necessary supplement to the classic ALM (\ref{ALM-E}), especially for the case where $\hbox{Prox}_\theta^r(x)$ has a closed-form representation but $\|A^TA\|$ is large.

The rest of this paper is organized as follows. We state the model to be considered and generalize the balanced ALM (\ref{SD-PPA}) for this model in Section \ref{Sec:algorithm}. Then, we conduct convergence analysis for the balanced ALM in Section \ref{Sec:convergence}. In Section \ref{Sec:splitting}, we extend our discussion to separable convex programming models and propose a splitting version of the balanced ALM. An alternative strategy for balancing is discussed in Section \ref{Sec:balance}. In Section \ref{Sec:remarks}, from the PPA perspective, we briefly discuss how to further generalize the algorithms to be proposed in Sections \ref{Sec:algorithm}-\ref{Sec:balance}. Finally, some conclusions are made in Section \ref{Sec:conclusion}.

\section{Model and algorithm}\label{Sec:algorithm}

\setcounter{equation}{0}

Note that the classic ALM (\ref{ALM-E}) was proposed in the context of the canonical convex programming model with linear equality constraints (\ref{Problem-E}). Despite that our initial aim is to consider the balanced ALM (\ref{SD-PPA}) for the model (\ref{Problem-E}), the balanced ALM (\ref{SD-PPA}) does can be generalized to the more general convex programming model with both linear equality and inequality constraints. We thus present our work in a more general setting as below.

\subsection{Model}

Instead of  \eqref{Problem-E}, let us consider the following more general convex programming model with both linear equality and inequality constraints:
\[  \label{Problem-B1}
    \min \{{\theta}(x) \; | \;  Ax =b\  (\hbox{or} \ge b), \, x\in {\cal X} \},
   \]
in which the setting is same as (\ref{Problem-E}) except that the linear inequality $Ax\ge b$ is also included. The solution set of (\ref{Problem-B1}) is assumed to be nonempty throughout our discussion. With the consideration of the more general model (\ref{Problem-B1}), the applicable range of the algorithm to be proposed is thus wider than that of the classic ALM (\ref{ALM-E}).

\subsection{Algorithm}

Now, let us generalize the balanced ALM (\ref{SD-PPA}) to the model (\ref{Problem-B1}), and the name remains for simplicity. Recall that our main purpose is to balance the subproblems in the classic ALM (\ref{ALM-E}).

{\hdp{
{\bf Algorithm: the balanced ALM for (\ref{Problem-B1})}
\smallskip

Let $r>0$ and $\delta>0$ be arbitrary constants; $H_0$ be defined in (\ref{H0}). Denote
$$
   q_0^k:=x^k + \frac{1}{r}A^T\lambda^k\;\;\hbox{and}\;\; s_0^k:=A(2{x}^{k+1}-x^k) - b.
$$
Then, with $(x^k, \lambda^k)$, the new iterate $(x^{k+1},{\lambda}^{k+1})$ is generated via the following steps:
\begin{subequations} \label{SD-PPA1-P}
\begin{numcases}{\hskip-1.5cm\;}
\label{SD-xP}   {x}^{k+1} =
  \arg\min  \bigl\{   \theta(x)  + \frac{r}{2}\|x- q_0^k\|^2 \; |\;  x\in {\cal X} \bigr\},  \\
\label{SD-yP}    {\lambda}^{k+1} =  \arg\min \bigl\{  \frac{1}{2}(\lambda-\lambda^k)H_0(\lambda-\lambda^k)+(s_0^k)^T\lambda \;|\; \lambda\in \Lambda\bigr\}.
\end{numcases}
\end{subequations}
}}

\begin{remark}\label{easyH}
It is easy to see that the balanced ALM (\ref{SD-PPA1-P}) is reduced to the aforementioned (\ref{SD-PPA}) if the model (\ref{Problem-E}) is considered. In particular, we have $\Lambda=\Re^m$ and thus the subproblem (\ref{SD-yP}) is reduced to finding $\lambda^{k+1}$ such that
$$
H_0(\lambda-\lambda^k)= -s_0^k.
$$
When $Ax\ge b$ is considered in \eqref{Problem-B1}, we have $\Lambda=\Re^m_+$. For this case, the subproblem \eqref{SD-yP} is reduced to the standard quadratic programming with non-negative sign constraints
$$
\min \bigl\{\frac{1}{2}(\lambda-\lambda^k)H_0(\lambda-\lambda^k)+(s_0^k)^T\lambda \;|\; \lambda \in \Re^n_+\bigl\},
$$
or equivalently, the linear complementarity problem
$$ 0\le \lambda \; \;  {\boldsymbol{\perp}} \;\;   \bigl\{ H_0 (\lambda - \lambda^k)   + s_0^k \} \ge 0.   $$
Recall that the matrix $H_0$ defined in (\ref{H0}) is positive definite and it can be well conditioned with appropriate choices of $r$ and $\delta$. Hence, it is extremely easy to decompose $H_0$, e.g., by the Cholesky decomposition. Then, many benchmark solvers including the well-known Lemke algorithm and conjugate gradient method, can be found in various textbooks (e.g., \cite{NW,GL,Stoer}), monographs (e.g., \cite{CPS}), and papers (e.g., \cite{He94MP,He94NM}).
\end{remark}

\begin{remark}

Recall that the balanced ALM (\ref{SD-PPA1-P}) is featured by the fact that the function $\theta(x)$ and the coefficient matrix $A$ are decoupled without any explicit or implicit condition related to $A$ in (\ref{SD-xP}). Compared with the PALM (\ref{PALM-E}) and the primal-dual method (\ref{CP-PPA}), the balanced ALM (\ref{SD-PPA1-P}) also has two parameters, $\delta$ and $r$, whose only restriction is their sign. As we will show in Section \ref{Sec:convergence}, the only essential role of $\delta$ is to theoretically ensure the positive definiteness of the corresponding matrix $H$ as defined in (\ref{PPA-T1H}). Therefore, there is no particular motivation to tune $\delta$ for different applications, and it can be just fixed as a small value beforehand. For the parameter $r$, just as the same parameter in the classic ALM (\ref{ALM-E}), there is full-extent flexibility to tune this parameter. Certainly, how to tune $r$ depends on the specific model and dataset under discussion, whilst there is no generic and unified theory to determine the optimal choice for all cases.
\end{remark}

\section{Convergence analysis}\label{Sec:convergence}

\setcounter{equation}{0}

In this section, we prove the convergence of the balanced ALM (\ref{SD-PPA1-P}), and estimate its worst-case convergence rate measured by the iteration complexity.

\subsection{Variational inequality characterization of (\ref{Problem-B1})}\label{Sec:convergence-VI}

Following our previous works \cite{HY-SIIMS,HY-SINUM}, our analysis will be conducted in the variational inequality (VI) context. We first derive the VI characterization for the optimality condition of the model (\ref{Problem-B1}). Let the Lagrangian function of the problem (\ref{Problem-B1}) be defined as
\[   \label{Lagrangian-EI}
      L(x,\lambda)= {\theta}(x) -  \lambda^T(Ax-b), \]
with $\lambda \in \Re^m$ the Lagrange multiplier. Since both linear equality and inequality constraints are considered in (\ref{Problem-B1}), let us define
  \[\label{OL2}   \Omega:= {\cal X}\times \Lambda \quad \hbox{where}   \quad     \Lambda :=\biggl\{ \begin{array}{ll}
               \Re^m,           &   \hbox{if   $Ax= b$} ,   \\
                \Re^m_+,    &        \hbox{if   $Ax \ge b$}.
                \end{array}
    \]
The  pair $(x^*,\lambda^*)\in \Omega$ is called a saddle point of the Lagrangian function (\ref{Lagrangian-EI}) if it satisfies the inequalities
 $$    L_{\lambda\in\Lambda}(x^*,\lambda) \le L(x^*,\lambda^*) \le L_{x\in {\cal X}}(x,\lambda^*). $$
Alternatively, we can write these inequalities as the following VIs:
\[   \label{VI1-Chara}
    \left\{ \begin{array}{lrl}
     x^*\in {\cal X}, & {\theta}(x) - {\theta}(x^*) + (x-x^*)^T(- A^T\lambda^*) \ge 0,
      & \forall\, x\in {\cal X},
      \\
  \lambda^*\in \Lambda,  &(\lambda -\lambda^*)^T(Ax^*-b)\ge 0,  &  \forall \; \lambda \in \Lambda,
        \end{array} \right.
       \]
or in the compact format
\begin{subequations}  \label{VI1}
\[  \label{VI1-F}
    w^*\in \Omega, \quad {\theta}(x) -{\theta}(x^*) + (w-w^*)^T F(w^*) \ge 0, \quad \forall \,  u\in
      \Omega,
     \]
where
\[ \label{VI1-CM}
     w = \left(\begin{array}{c}
                     x\\
                  \lambda  \end{array} \right),
  \quad
    F(w) =\left(\begin{array}{c}
     - A^T\lambda  \\
     Ax-b \end{array} \right) \quad
  \hbox{and} \quad      \Omega= {\cal X} \times \Lambda.
   \]
   \end{subequations}
 Note that the operator $F$ defined in \eqref{VI1-CM} is affine with a skew-symmetric matrix and thus we have
\[ \label{Skew-S0} (w-\tilde{w})^T(F(w)-F(\tilde{w}))\equiv 0.  \]
We also call \eqref{VI1} a monotone mixed variational inequality because the function ${\theta}$ is convex and the operator $F$ has the property
 \eqref{Skew-S0}. We denote by $\Omega^*$ the solution set of the VI \eqref{VI1}; it is also the set of the saddle points of the Lagrangian function \eqref{Lagrangian-EI}.

\subsection{Contraction}

We need to show that the sequence generated by the balanced ALM (\ref{SD-PPA1-P}) is contractive with respect to ${\Omega^*}$, the solution set of the VI (\ref{VI1}). This is the key property to ensure its convergence. Before that, let us recall a basic lemma whose proof is elementary and can be found in, e.g., \cite{Beck}.

\begin{lemma} \label{CP-TF}
\begin{subequations} \label{CP-TF0}
 Let ${\cal X}\subset \Re^n$ be a closed convex set, $\theta(x)$ and $f(x)$ be convex functions.
   If $f$ is differentiable, and the solution set of the minimization problem
   $$
   \min\{\theta(x) + f(x) \, |\,x\in {\cal X}\}
   $$ is nonempty, then it holds that
  \[  \label{CP-TF1}   x^*  \in \arg\min \{ \theta(x) + f(x) \, |\, x\in {\cal X}\}
     \]
if and only if
\[\label{CP-TF2}
      x^*\in {\cal X}, \quad   \theta(x) - \theta(x^*) + (x-x^*)^T\nabla f(x^*) \ge 0, \quad \forall\, x\in {\cal X}.
      \]
\end{subequations}
\end{lemma}

To show the contraction property of the sequence generated by the balanced ALM (\ref{SD-PPA1-P}), the first step is to fathom the difference of an iterate generate by the balanced ALM (\ref{SD-PPA1-P}) from a solution point $w^* \in \Omega^*$. Recall the definition of $H_0$ in (\ref{H0}). Let us define
\[ \label{PPA-T1H}
  H  = \left(\begin{array}{cc}
          r I_n      &      A^T \\
            A                &    \displaystyle{\frac{1}{r}}AA^T + \delta I_m
            \end{array}\right)=\left(\begin{array}{cc}
          r I_n      &      A^T \\
            A                &   H_0
            \end{array}\right).
            \]

\begin{proposition}\label{PropH}
The matrix $H$ defined in (\ref{PPA-T1H}) is positive definite.
\end{proposition}
\noindent{\bf Proof}.
Notice that
$$    H=    \left(\begin{array}{cc}
          r I_n      &      A^T \\
            A                &    \displaystyle{\frac{1}{r}}AA^T
            \end{array}\right)   +   \left(\begin{array}{cc}
             0    &     0 \\
             0             &   \delta I_m
            \end{array}\right)  =  \left(\begin{array}{c}
         \sqrt{r}  I_n   \\[0.1cm]
              \sqrt{\frac{1}{r}}A
            \end{array}\right) \Bigl( \sqrt{r} I_n,  \textstyle{\sqrt{\frac{1}{r}}}A^T \Bigr)+ \left(\begin{array}{cc}
             0    &     0 \\
             0             &   \delta I_m
            \end{array}\right),  $$
   for any $w=(x,\lambda)\ne 0$. Thus, we have
   $$    w^THw=    \big\|\sqrt{r}x + \textstyle{ \sqrt{\frac{1}{r}}}A^T\lambda\big\|^2  + \delta \|\lambda\|^2>0,  $$
       and therefore the matrix $H$ is positive definite.  \hfill {$\Box$}

In the following theorem, we will express the difference of an iterate generated by the balanced ALM (\ref{SD-PPA1-P}) from a solution point $w^* \in \Omega^*$ in the context of VIs.

\begin{theorem}  \label{HeB1} Let $\{w^k=(x^k,\lambda^k)\}$ be the sequence generated by the balanced ALM \eqref{SD-PPA1-P} and $H$ be defined in (\ref{PPA-T1H}). Then we have
\[  \label{PPA-T1F}  {w}^{k+1} \in \Omega, \;\;
      \theta(x)
      -\theta({x}^{k+1})  +(w-  {w}^{k+1})^TF({w}^{k+1}) \ge  (w- {w}^{k+1})^T H(w^k-{w}^{k+1}) , \quad \forall\, w\in
         \Omega.
          \]
\end{theorem}
\noindent{\bf Proof}.
According to Lemma \ref{CP-TF}, the solution ${x}^{k+1}$ of the subproblem  \eqref{SD-xP}  can be characterized by the VI
$$ {x}^{k+1}\in {\cal X},  \;\;
   {\theta}(x) - {\theta}({x}^{k+1}) +
 (x-{x}^{k+1})^T\bigl\{-A^T\lambda^k  +   r ({x}^{k+1} - x^k)  \bigr\}\ge 0,  \;\; \forall\, x\in {\cal X}.
 $$
Then, for any unknown ${\lambda}^{k+1}$, we have
\begin{eqnarray} \label{SD-Pre-x}
 \lefteqn{ {x}^{k+1}\in {\cal X},  \;\;
   {\theta}(x) - {\theta}({x}^{k+1}) +
 (x-{x}^{k+1})^T(-A^T{\lambda}^{k+1} )  } \nn \\
   & & \qquad \qquad\qquad  \ge  (x-{x}^{k+1})^T  \bigl\{ r (x^k - {x}^{k+1})  + A^T(\lambda^k - {\lambda}^{k+1}) \bigr\},  \;\; \forall\, x\in {\cal X}.
   \end{eqnarray}
Similarly, because of Lemma \ref{CP-TF},  the solution ${\lambda}^{k+1}$ of the subproblem \eqref{SD-yP}  can be characterized by the VI
  $$   {\lambda}^{k+1}\in \Lambda ,  \;\;
  (\lambda-{\lambda}^{k+1})^T\Bigl\{ \Bigl( A[2{x}^{k+1}- x^k] -b\Bigr)  +  H_0 ({\lambda}^{k+1}-\lambda^k) \Bigr\}
     \ge 0, \quad \forall\, \lambda\in \Lambda.  $$
Recall the definition of $H_0$ in (\ref{H0}). We thus have
\begin{eqnarray} \label{SD-Pre-l}
 \lefteqn{ {\lambda}^{k+1}\in \Lambda,  \;\;
 (\lambda-{\lambda}^{k+1})^T(A{x}^{k+1} -b)  } \nn \\
   & & \qquad \quad  \ge  (\lambda-{\lambda}^{k+1})^T  \bigl\{ (A(x^k - {x}^{k+1})  + \Bigl(\displaystyle{\frac{1}{r}}AA^T + \delta I_m \Bigr) (\lambda^k - {\lambda}^{k+1}) \bigr\},  \;\; \forall\, \lambda\in \Lambda .
   \end{eqnarray}
Combining  \eqref{SD-Pre-x} and  \eqref{SD-Pre-l}, and  using the notation in \eqref{VI1}, we obtain the assertion (\ref{PPA-T1F}).
  \hfill {$\Box$}

In the following theorem, we will prove an important inequality which measures the difference of an iterate generated by the balanced ALM (\ref{SD-PPA1-P}) from a solution point $w^* \in \Omega^*$ more explicitly by $H$-norm-induced distances. This inequality is also the basis of estimating the convergence rate measured by the iteration complexity for the balanced ALM  (\ref{SD-PPA1-P}).

\begin{theorem} \label{HB1-A}
Let $\{w^k=(x^k,\lambda^k)\}$ be the sequence generated by the balanced ALM \eqref{SD-PPA1-P} and $H$ be defined in (\ref{PPA-T1H}). Then we have
\begin{eqnarray} \label{Hau1-A0}
\lefteqn{  {\theta}(x) - {\theta}({x}^{k+1}) +(w- {w}^{k+1})^T
             F(w)   }  \nn \\
             &  \ge  &  \frac{1}{2}
    \bigl(\|w-w^{k+1}\|_{H}^2 -\|w-w^k\|_{H}^2  \bigr) +
        \frac{1}{2}\|w^k-{w}^{k+1}\|_{H}^2, \;\;  \forall w\in \Omega.
        \end{eqnarray}
    \end{theorem}
\noindent{\bf Proof}. It follows from \eqref{Skew-S0} that
$$
(w- {w}^{k+1})^T  F({w}^{k+1}) =(w- {w}^{k+1})^T  F(w),
$$
and thus the  left-hand side of \eqref{PPA-T1F} equals
 $$  {\theta}(x) - {\theta}({x}^{k+1}) +(w- {w}^{k+1})^T F(w). $$
Consequently, because of \eqref{PPA-T1F},  we get
\begin{equation} \label{LEM-M-H}
   {w}^{k+1} \in \Omega, \quad   {\theta}(x) - {\theta}(\tilde{x}^k) + (w- {w}^{k+1} )^T F(w)
      \ge (w- {w}^{k+1})^T H(w^k -
    w^{k+1}),  \quad  \forall w\in \Omega.
   \end{equation}
Applying the identity
$$
      b^TH(b-a) = \frac{1}{2} \{\|b\|_{H}^2 -\|a\|_{H}^2\} +\frac{1}{2}\|a-b\|_{H}^2
$$
to the right-hand side of \eqref{LEM-M-H} with $  a=w-w^k$ and $ b=w-{w}^{k+1}$, we thus obtain
\[ \label{LEM-M1-2}
(w- {w}^{k+1})^T H(w^k - w^{k+1})
     \, =\,
        \frac{1}{2} \bigl(
   \|w-w^{k+1}\|_{H}^2 - \|w-w^k\|_{H}^2 \bigr)
  +  \frac{1}{2}\|w^k-{w}^{k+1}\|_{H}^2 .
   \]
Substituting \eqref{LEM-M1-2}  into the right-hand side of \eqref{LEM-M-H}, we prove the assertion (\ref{Hau1-A0}).  \hfill {$\Box$}

Now, with Theorems \ref{HeB1} and \ref{HB1-A}, the contraction property of the sequence generated by the balanced ALM (\ref{SD-PPA1-P}) with respect to $\Omega^*$ can be proved.

\begin{theorem} \label{HB1-B} Let $\{w^k=(x^k,\lambda^k)\}$ be the sequence generated by the balanced ALM \eqref{SD-PPA1-P} and $H$
be defined in (\ref{PPA-T1H}). Then we have
\[  \label{Hau1-Bw}   \|w^{k+1}-w^*\|_H^2  \le \|w^k-w^* \|_H^2  -\|w^k-{w}^{k+1}\|_{H}^2, \;\;
 \forall w^*\in \Omega^*.
\]
\end{theorem}
\noindent{\bf Proof}.   Setting  $w$ in \eqref{Hau1-A0}  as any fixed $w^*\in \Omega^*$, we get
\begin{eqnarray*}
   \lefteqn{ \|w^k-w^*\|_H^2  -\|w^{k+1}-w^*\|_{H}^2  - \|w^k-{w}^{k+1}\|_{H}^2 }\nn \\
      &\ge&  2\bigl\{ {\theta}({x}^{k+1})-\theta(x^*) +({w}^{k+1}-w^*)^T
             F(w^*) \bigr\}, \;\;   \forall w^*\in \Omega^*.
        \end{eqnarray*}
Since $w^*\in \Omega^*$ and ${w}^{k+1} \in \Omega$,   according to \eqref{VI1},  the right-hand side of the
last inequality is non-negative. Thus,  the assertion of this theorem follows directly. \hfill {$\Box$}

\medskip

\subsection{Convergence}

With the contraction property established in Theorem \ref{HB1-B}, it is easy to prove the convergence of the sequence $\{w^k\}$ generated by the balanced ALM \eqref{SD-PPA1-P}.

 \begin{theorem} \label{HauptB} Let $\{w^k=(x^k,\lambda^k)\}$ be the sequence generated by the balanced ALM \eqref{SD-PPA1-P} and $H$ be defined in (\ref{PPA-T1H}).  Then, the sequence $\{w^k\}$ converges to  some
  $w^{\infty} \in \Omega^*$.
    \end{theorem}
 \noindent{\bf Proof}.  First of all,  it follows from \eqref{Hau1-Bw}    that the sequence  $\{w^k\}$
   is bounded and
   \[   \label{Key-Eq} \lim_{k\to\infty}  \| w^k - {w}^{k+1}  \|_H^2 =0.
              \]
  Let $w^{\infty}$ be a  cluster point of $\{w^k\}$ and
 $\{w^{k_j}\}$ be a subsequence converging to $w^{\infty}$. It follows from  \eqref{PPA-T1F}  that
 $$    {w}^{k_j}\in \Omega,  \;\; \theta(x) - \theta({x}^{k_j})  +  (w- {w}^{k_j})^T F({w}^{k_j}) \ge
      (w  -{w}^{k_j})^TH(w^{k_j-1}-{w}^{k_j}), \quad \forall  \, w\in {\Omega}. $$
 Since the matrix $H$ is positive definite, it follows from \eqref{Key-Eq} and the continuity of $\theta(x)$ and $F(w)$ that
 $$  w^{\infty}\in \Omega,  \;\; \theta(x) - \theta(x^{\infty})  +  (w- w^{\infty})^T F(w^{\infty}) \ge
    0, \quad \forall  \, w \in {\Omega}.   $$
This VI above indicates that $w^{\infty}$ is a solution point of \eqref{VI1}.
 Finally, because of \eqref{Hau1-Bw}, we have
 $$     \|w^{k+1} - w^{\infty} \|_H^2 \le \|w^k - w^{\infty}\|_H^2,   $$
 and thus $\{w^k\}$ converges to $w^{\infty}$.  The proof is complete. \hfill {$\Box$}

 \subsection{Convergence rate}

Following the VI-based technique established in our earlier work \cite{HY-SINUM}, we can estimate the worst-case $O(1/t)$ convergence rate measured by the iteration complexity for the balanced ALM \eqref{SD-PPA1-P} where $t$ is the iteration counter.

Let us recall some necessary details which can also be found in \cite{HY-SINUM}. If $\tilde{w}$ is a solution point of the VI \eqref{VI1}, then we have
$$
   \tilde{w}\in \Omega, \quad  \theta (x) - \theta(\tilde{x}) + (w-\tilde{w})^T F(\tilde{w}) \ge 0, \quad  \forall \, w\in \Omega.
  $$
Because of \eqref{Skew-S0}, $\tilde{w}$ also satisfies
$$
 \tilde{w}\in \Omega, \quad    \theta (x) - \theta(\tilde{x}) + (w-\tilde{w})^T F(w) \ge 0, \quad \forall \, w\in \Omega.
  $$
Thus, for given $\epsilon>0$,
  $\tilde{w} \in {\Omega}$ is called  an $\epsilon$-approximate solution of VI \eqref{VI1} if it satisfies
\[   \label{Aproxi-E}
 \tilde{w}\in \Omega, \;\; \theta(x) - \theta(\tilde{x}) + (w-\tilde{w})^TF(w)\ge -\epsilon,
\;\;\forall \; w \in {\cal D}_{(\tilde{w})},
   \]
where
$$ {\cal D}_{(\tilde{w})} =\{  w\in \Omega \, |\, \|w-\tilde{w}\|\le 1 \}.
   $$
Thus, to establish the worst-case $O(1/t)$ convergence rate for the balanced ALM \eqref{SD-PPA1-P}, we need to show that, for given $\epsilon>0$, after $t$ iterations, we can find $\tilde{w} \in \Omega$, such that
\[ \label{Main-ATT} \tilde{w}\in \Omega,\quad \hbox{and} \quad  \sup_{w\in
 {\cal D}_{(\tilde{w})}} \bigl\{
     \theta(\tilde{x}) - \theta(x) +
        (\tilde{w}-w)^T F(w) \bigr\}\le \epsilon=O(1/t).
 \]
We present this result in the following theorem.

\begin{theorem} \label{HauptS}
Let $\{w^k=(x^k,\lambda^k)\}$ be the sequence generated by the balanced ALM \eqref{SD-PPA1-P} and $H$ be defined in (\ref{PPA-T1H}). For any integer number $t>0$, if we define
 \begin{equation} \label{Ergoticw}
\tilde{w}_t := \frac{1}{t+1} \sum_{k=0}^t{w}^{k+1},
\end{equation}
then we have
 \begin{equation} \label{THMXS0}
 \tilde{w}_t\in\Omega, \quad \theta(\tilde{x}_t) -\theta(x)  + (\tilde{w}_{t}-w)^T F(w) \le
\frac{1}{2(t+1)}\|w-w^0\|_H^2,\quad  \forall w\in\Omega.
\end{equation}
\end{theorem}
\noindent{\bf Proof}. First,  it follows from  \eqref{Hau1-A0} that,  for all $k\ge0$, we have
\begin{equation}
\label{HH-SE}  {w}^{k+1} \in \Omega, \quad   \theta(x)\! -\! \theta({x}^{k+1} )\! + \! (w- {w}^{k+1} )^TF(w)
       + \frac{1}{2}\|w-w^k\|_H^2 \ge \frac{1}{2}\|w-w^{k+1}\|_H^2, \; \forall w\in\Omega.
    \end{equation}
Summarizing the inequalities \eqref{HH-SE} over $k=0,1,\ldots, t$, we
obtain
$$  (t+1)\theta(x)
 - \sum_{k=0}^t \theta({x}^{k+1}) + \Bigl(   (t+1)  w
 - \sum_{k=0}^t {w}^{k+1}\Bigr)^T
    F(w) +  \frac{1}{2}\|w-w^0\|_H^2 \ge 0, \quad \forall w\in\Omega.
   $$
It follows from (\ref{Ergoticw}) that 
\begin{equation} \label{THMXS0-wA}
  \frac{1}{t+1} \sum_{k=0}^t \theta({x}^{k+1}) -\theta(x)  + (\tilde{w}_{t}-w)^T
F(w) \le \frac{1}{2(t+1)}\|w-w^0\|_H^2,\quad  \forall w\in\Omega.
\end{equation}
Note that $\tilde{w}_{t}$ defined in (\ref{Ergoticw}) is a convex combination of all iterates $w^k$ for $k=0,\cdots,t$, and $\theta(x)$ is convex. We thus have
 $$  \tilde{x}_{t}= \frac{1}{t+1} \sum_{k=0}^t{x}^{k+1},
  $$
and also
$$ \theta(\tilde{x}_t) \le\frac{1}{t+1}\sum_{k=0}^t
\theta({x}^{k+1}).
  $$
Substituting it into \eqref{THMXS0-wA}, the assertion (\ref{THMXS0}) of this theorem follows directly. \hfill \quad {\normalsize$\Box$}

Then, because of (\ref{Main-ATT}), the inequality (\ref{THMXS0}) indicates that ${\tilde w}_t$ defined in (\ref{Ergoticw}), which is the average of the first $t$ iterates generated by the balanced ALM \eqref{SD-PPA1-P}, is an approximate solution of the VI (\ref{VI1}) with an accuracy of $O(1/t)$. Hence, the worst-case $O(1/t)$ convergence rate measured by the iteration complexity is established for the balanced ALM \eqref{SD-PPA1-P} in the ergodic sense.

\section{Splitting versions of the balanced ALM (\ref{SD-PPA1-P}) for separable convex programming}\label{Sec:splitting}

\setcounter{equation}{0}

The classic ALM (\ref{ALM-E}) plays an extremely influential role in solving various separable cases of the generic model (\ref{Problem-E}) when the objective function of such a model can be represented as the sum of multiple functions without coupled variables. For these separable models, the classic ALM (\ref{ALM-E}) has been adapted into various splitting versions by decomposing the primeval $x$-subproblem (\ref{ALM-E-x}) into smaller ones. These splitting versions take advantage of the separable structure in the model more effectively; the decomposed subproblems are usually easier in the sense that
each of them only needs to tackle one function component. For various applications including the mentioned sparsity- and low-rank-promoted ones in data science domains, splitting versions of the ALM (\ref{ALM-E}) may generate subproblems that are easy enough to have closed-form solutions. Among various splitting versions of the classic ALM (\ref{ALM-E}), the most popular one is probably the mentioned ADMM in \cite{GM75}, which suggests splitting the $x$-subproblem (\ref{ALM-E-x}) into two sequentially when the model (\ref{Problem-E}) has a two-block separable structure.

In this section, in parallel with the successful legacy of the classic ALM (\ref{ALM-E}) and its splitting versions, we also discuss how to design splitting versions for the balanced ALM \eqref{SD-PPA1-P} when the model (\ref{Problem-B1}) is separable. For succinctness and without ambiguity, we reuse some letters and notation as those in Sections \ref{Sec:algorithm} and \ref{Sec:convergence}.

\subsection{Model}

Let us consider the separable convex programming model with both linear equality and inequality constraints
\[  \label{Problem-m}
  \min \Bigl\{ \sum_{i=1}^{p} \theta_i(x_i)   \;\big|\;   \sum_{i=1}^{p} A_ix_i=b\ (\hbox{or} \ge b) ,  \;\;  x_i\in {\cal X}_i \Bigr\},\]
where $\theta_i: {\Re}^{n_i}\to {\Re}, \, i=1,\ldots, p$, are closed
proper convex functions and they are not necessarily smooth; ${\cal
X}_i\subseteq \Re^{n_i},\, i=1,\ldots, p$, are closed convex sets;
$A_i\in \Re^{m\times n_i},\, i=1,\ldots, p$, are given matrices; and $b\in \Re^m$ is a given vector. The model (\ref{Problem-m}) can be regarded as an extension of the model (\ref{Problem-B1}) from $p=1$ to $p\ge 1$. Let us only consider the multiple-block separable case with $p\ge2$ for (\ref{Problem-m}). Similarly as (\ref{OL2}), we reuse the letters and define
\[\label{OL-p}
     \Omega = \prod_{i=1}^p {\cal X}_i \times \Lambda \qquad \hbox{where}\qquad   \Lambda  =\left\{ \begin{array}{ll}
               \Re^m,           &   \hbox{if   $\sum_{i=1}^{p} A_ix_i =  b$} ,   \\[0.2cm]
                \Re^m_+,    &        \hbox{if   $\sum_{i=1}^{p} A_ix_i\ge b$}.
                \end{array}        \right.
                \]

\subsection{Algorithm}

Now, we extend the balanced ALM (\ref{SD-PPA1-P}) to the multiple-block separable convex programming model (\ref{Problem-m}) and present a splitting version of (\ref{SD-PPA1-P}) below.

{\hdp{
{\bf Algorithm: A splitting version of the balanced ALM (\ref{SD-PPA1-P}) for (\ref{Problem-m})}
\smallskip

Let $r_i>0$ for $i=1,2,\cdots, p$, and $\delta>0$ be arbitrary constants. Define
\[\label{Hp}
H_p=\sum_{i=1}^p \frac{1}{r_i}A_iA_i^T + \delta I_m,
\]
$$
q^k_i:=x_i^k + \frac{1}{r_i}A^T\lambda^k, \; \hbox{for}\; i=1,2,\cdots,p;
\;\hbox{and}\;
s^k=\sum_{i=1}^pA_i(2x_i^{k+1}-x_i^k)-b.
$$
Then, with $w^k= (x_1^{k},x_2^{k},\cdots, x_p^{k}, \lambda^{k})$, the new iterate ${w}^{k+1} =({x}_1^{k+1}, {x}_2^{k+1},\cdots, {x}_p^{k+1}, {\lambda}^{k+1})$ is generated via the following steps:

\begin{subequations} \label{M-PRE}
\begin{numcases}{\hskip-1.5cm\;}
\label{M-SD-x}   {x}_i^{k+1}\in\arg\min\bigl\{\theta_i(x_i)  +\frac{r_i}{2} \|x_i-q_i^k\|^2 \; | \;  x_i\in {\cal X}_i \bigr\}, i=1,2,\cdots, p;\\
\label{M-SD-y}   {\lambda}^{k+1} =  \arg\min \left\{ {\displaystyle\frac{1}{2}}(\lambda-\lambda^k)^TH_p(\lambda-\lambda^k)                                                                                         + (s^k)^T\lambda \; | \; \lambda\in \Lambda\right\}.
\end{numcases}
\end{subequations}
}}

\begin{remark}
The subproblems in (\ref{M-PRE}) are of the same structure as those in (\ref{SD-PPA1-P}). For the $x_i$-subproblem (\ref{M-SD-x}), the function $\theta_i(x_i)$ and the coefficient $A_i$ are decoupled without any explicit or implicit condition related to $A_i$, and thus it is also reduced to estimating the proximity operator of $\theta_i(x_i)$ when ${\cal X}_i=\Re^{n_i}$. In addition, the $\lambda$-subproblem (\ref{M-SD-y}) is a positive definite system of linear equations or a standard quadratic programming with non-negative sign constraints. Note that all $r_i$'s have no other restriction than the sign requirement. Hence, the algorithm (\ref{M-PRE}) keeps all features of the balanced ALM (\ref{SD-PPA1-P}).
\end{remark}

\begin{remark}
For generality, we consider different $r_i$ for different $x_i$-subproblems. They can be identical for simplicity. Similarly as the balanced ALM (\ref{SD-PPA1-P}), to implement the algorithm (\ref{M-PRE}), empirically we can fix $\delta$ as a small positive constant throughout.
\end{remark}

\subsection{Convergence analysis}

In this subsection, we follow the analysis in Section \ref{Sec:convergence} and prove the convergence of the splitting version of the balanced ALM (\ref{M-PRE}).

\subsubsection{Variational inequality characterization of (\ref{Problem-m})}

For convergence analysis purpose, we also need the VI characterization for the optimality condition of the model (\ref{Problem-m}).
Let $\lambda\in\Re^m$ be the Lagrange multiplier of (\ref{Problem-m}) and the Lagrangian  function  of the problem \eqref{Problem-m} be defined as
\[ \label{Lagrange-F}
  L(x_1,\ldots,x_p,\lambda) =  \sum_{i=1}^{p}\theta_i(x_i) -\lambda^T\Bigl(\sum_{i=1}^{p} A_ix_i-b\Bigr).\]
Similarly as Section \ref{Sec:convergence-VI}, we reuse the letters and know that finding a saddle point of $ L(x_1,\ldots,x_p,\lambda)$ can be written as the following VI:
\begin{subequations} \label{VIp-FORM}
\[  \label{VI-FORM-Q}
    w^*\in\Omega, \quad \theta(x) - \theta(x^*) + (w-w^*)^T F(w^*) \ge0, \quad
\forall\, w\in \Omega,
    \]
where
\[  \label{VI-FORM-F}    w=\left(\!\!\begin{array}{c}
             x_1\\
         \vdots \\
            x_{p} \\[0.1cm]
            \lambda
             \end{array}\!\! \right), \quad
   x=\left(\!\!\begin{array}{c}
             x_1\\
         \vdots \\
            x_{p} \\
             \end{array}\!\! \right),
    \quad   \theta(x) = \sum_{i=1}^{p} \theta_i(x_i), \quad
       F(w) = \left(\!\!\begin{array}{c}
    - A_1^T\lambda\\
       \vdots \\
    -A_{p}^T\lambda \\[0.1cm]
       \sum_{i=1}^{p} A_ix_i-b
    \end{array}\!\! \right),
  \]
  and $\Omega$ is defined in (\ref{OL-p}).
\end{subequations}
Again, we denote by $\Omega^*$ the solution set of the VI \eqref{VIp-FORM}.

\subsubsection{Convergence}

Let us recall the  proofs in Section \ref{Sec:convergence} for the convergence of the balanced ALM (\ref{SD-PPA1-P}). It is easy to see that the crucial step is to identify the difference between an iterate and a solution point by the inequality (\ref{PPA-T1F}) in Theorem \ref{HeB1}, in which the matrix $H$ should be positive definite as proved in Proposition \ref{PropH} so that the difference can be measured by distances defined by the $H$-norm. After Proposition \ref{PropH} and Theorem \ref{HeB1} are proved, the remaining part of the proof for the convergence as well as the worst-case convergence rate is subroutine. Hence, to prove the convergence of the splitting version of the balanced ALM (\ref{M-PRE}), we only need to prove an inequality similar as (\ref{PPA-T1F}) in which the accompanying matrix is also positive definite.

\begin{proposition}\label{PropHp}
Let $r_i>0$ for $i=1,2,\cdots, p$, and $\delta>0$ be arbitrary constants. The matrix defined as
\[ \label{PPA-TpH}
  {H}  = \left(\begin{array}{ccccc}
         r_1I_{n_1}      &       0      & \cdots  &  0    &     A_1^T \\[0.2cm]
             0                                   & \quad \ddots\quad    & \ddots      & \vdots    & \vdots     \\[0.3cm]
         \vdots                              &   \ddots           & \quad \ddots \quad   &     0      & \vdots  \\[0.3cm]
              0                                  &     \cdots                          &       \quad  0 \quad       & r_pI_{n_p}  &   A_p^T \\
             A_1      & \cdots    &  \cdots     &   A_p             &    \displaystyle{\sum_{i=1}^p\frac{1}{r_i}}A_iA_i^T + \delta I_m
            \end{array}\right) \]
is positive definite.
\end{proposition}
\noindent{\bf Proof}. Note that
$$     H=     \sum_{i=1}^p H_i   +   \left(\begin{array}{cc}
           0     &     0 \\
             0             &   \delta I_m
            \end{array}\right) ,$$
where
$$  {H}_i  = \left(\begin{array}{ccccc}
                \quad        &  \quad     &         &   &    \\
                                         &  r_i I_{n_i}     &   & &  A_i^T      \\
                                        &                       &     &     \\           &                       &     &     \\
                   &         A_i            &   &  &   \displaystyle{\frac{1}{r_i}}A_iA_i^T
            \end{array}\right)  =  \left(\begin{array}{c}
                    \vdots   \\
                 \sqrt{r_i} I_{n_i}      \\[0.22cm]
                          \vdots                 \\[0.22cm]
                \displaystyle{\sqrt{\frac{1}{r_i}}}A_i
            \end{array}\right) \left(\begin{array}{ccccc}
                          \cdots         &      \sqrt{r_i} I_{n_i} &\; \cdots \;   &   \displaystyle{\sqrt{\frac{1}{r_i}}}A_i^T
            \end{array}\right).
  $$
  For any $ w=(x_1, \ldots, x_p, \lambda)\ne 0$, we have
            $$  w^THw=
               \sum_{i=1}^p \Big\|\sqrt{r_i}x_i + \sqrt{\frac{1}{r_i}}A_i^T\lambda\Big\|^2  + \delta \|\lambda\|^2>0.  $$
Hence, the matrix $H$  is positive definite. \hfill \quad {\normalsize$\Box$}

\begin{theorem}  \label{HeBp}  Let $\{{w}^{k}=(x_1^k,\cdots,x_p^k,\lambda^k)\}$ be the sequence generated by the balanced ALM \eqref{M-PRE} and $H$ be defined in (\ref{PPA-TpH}). Then, we have
\[  \label{PPA-wT}
{w}^{k+1} \in \Omega, \;\;
      \theta (x)
      -\theta({x}^{k+1})  +(w-  {w}^{k+1})^TF({w}^{k+1}) \ge  (w- {w}^{k+1})^TH(w^k- {w}^{k+1}) , \quad \forall\, w\in
         \Omega.
          \]
                 \end{theorem}
\noindent{\bf Proof}.  According to Lemma \ref{CP-TF}, for $i=1, 2, \ldots, p$, we have
$$ x_i^{k+1}\in {\cal X}_i, \quad
 \theta_i(x_i) - \theta_i(\tilde{x}_i^k)  + (x - {x}_i^{k+1})^T\{-A_i^T\lambda^k  + r_i ({x}_i^{k+1}-x_i^k)\} \ge 0,
    \;\; \forall \, x_i\in {\cal X}_i .   $$
 Then, for any unknown ${\lambda}^{k+1}$, we have
\begin{eqnarray} \label{M-SD-Pre-x}
 \lefteqn{ {x}_i^{k+1}\in {\cal X}_i,  \;\;
   {\theta_i}(x_i) - {\theta_i}({x}_i^{k+1}) +
 (x_i-{x}_i^{k+1})^T(-A^T{\lambda}^{k+1} )  } \nn \\
   & & \qquad \qquad\qquad  \ge  (x_i-{x}_i^{k+1})^T  \bigl\{ r_i (x_i^k - {x}_i^{k+1})  + A^T(\lambda^k - {\lambda}^{k+1}) \bigr\},
   \;\; \forall\, x_i\in {\cal X}_i.
   \end{eqnarray}
Also because of Lemma \ref{CP-TF},  ${\lambda}^{k+1}$ generated by \eqref{M-SD-y} is characterized by the VI
  $$   {\lambda}^{k+1}\in \Lambda ,  \;\;
  (\lambda-{\lambda}^{k+1})^T\Bigl\{ \Bigl( \sum_{i=1}^p A_i[2{x}_i^{k+1}- x_i^k] -b\Bigr)  + \Bigl(\sum_{i=1}^p\frac{1}{r_i}A_iA_i^T + \delta I_m \Bigr) ({\lambda}^{k+1}-\lambda^k) \Bigr\}
     \ge 0, \quad \forall\, \lambda\in \Lambda.  $$
It can be rewritten as
\begin{eqnarray} \label{M-SD-Pre-y}
 \lefteqn{ {\lambda}^{k+1}\in \Lambda,  \;\;
 (\lambda-{\lambda}^{k+1})^T\Bigl(\sum_{i=1}^p A_i{x}_i^{k+1} -b \Bigr)  }\nn \\
   & & \quad \ge  (\lambda-{\lambda}^{k+1})^T  \Bigl\{ \sum_{i=1}^p A_i(x_i^k - {x}_i^{k+1})  + \Bigl(\displaystyle{
      \sum_{i=1}^p\frac{1}{r_i}}A_iA_i^T+ \delta I_m \Bigr) (\lambda^k - {\lambda}^{k+1}) \Bigr\},  \forall\, \lambda\in \Lambda.
   \end{eqnarray}
Combining \eqref{M-SD-Pre-x} and  \eqref{M-SD-Pre-y}, and using the notation in \eqref{VIp-FORM}, we prove the assertion (\ref{PPA-wT}). \hfill \quad {\normalsize$\Box$}

As mentioned, based on Proposition \ref{PropHp} and Theorem \ref{HeBp}, similar conclusions as Theorems \ref{HB1-A}-\ref{HauptS} can be trivially proved. Thus, convergence results similar as those in Section \ref{Sec:convergence} can be obtained for the splitting version of the balanced ALM (\ref{M-PRE}); we omit the details for succinctness.

\section{An alternative strategy for balancing}\label{Sec:balance}
\setcounter{equation}{0}

The balanced ALM (\ref{SD-PPA1-P}) can be generalized to the splitting version (\ref{M-PRE}) if the model under discussion is changed from the one-block case (\ref{Problem-B1}) to the multiple-block case (\ref{Problem-m}). There are other ways for the generalization, in addition to the technique introduced in Section \ref{Sec:splitting}. In (\ref{M-PRE}), we see that each of the $x_i$-subproblems does not involve any quadratic term in form of $\frac{r_i}{2}\|A_ix_i-q_i^k\|^2$ so that it can be reduced to estimating the proximity operator of $\theta_i(x_i)$ when ${\cal X}_i=\Re^{n_i}$. In this sense, all such $x_i$-subproblems are preferred when it is easy to estimate the proximity operator of $\theta_i(x_i)$. On the other hand, all $A_i$'s are aggregated in the $\lambda$-subproblem (\ref{M-SD-y}) because of the matrix $H_p$ defined in (\ref{Hp}). For some cases where some or all $\|A_i^TA_i\|$ are large (or, some or all  $A_i$'s are ill-conditioned), it is preferred to consider alleviating the quadratic programming problem (\ref{M-SD-y}) by removing such $A_i^TA_i$ from $H_p$. Hence, from methodological point of view, it is also interesting to ask if we can keep terms in form of $\frac{r_i}{2}\|A_ix_i-q_i^k\|^2$ for some $x_i$-subproblems (\ref{M-SD-x}), and meanwhile remove the corresponding $A_iA_i^T$ from the matrix $H_p$ in (\ref{Hp}) so that the $\lambda$-subproblem (\ref{M-SD-y}) becomes easier. Accordingly, we propose to revise the splitting version of the balanced ALM (\ref{M-PRE}) such that some $x_i$-subproblems are in form of
$$
{x}_i^{k+1}\in\arg\min\bigl\{\theta_i(x_i)  +\frac{r_i}{2} \|A_ix_i-q_i^k\|^2 \; | \; x_i\in {\cal X}_i\big\},
$$
with $q_i^k$ a certain constant vector, and the corresponding $A_iA_i^T$ is excluded in the $\lambda$-subproblem (\ref{M-SD-y}).
This idea provides an alternative strategy for balancing the generated subproblems, and it enables a user to determine how to balance the difficulty of subproblems in accordance with the specific functions $\theta_i$'s, coefficient matrices $A_i$'s, and sets ${\cal X}_i$'s for a given application.

\subsection{Model}
For succinctness of notation, let us just take the special case of (\ref{Problem-m}) with $p=2$ and only linear equality constraints to illustrate our idea:
\[  \label{Problem-B2}
    \min \{\theta_1(x_1)  + \theta_2(x_2) \; | \;  A_1x_1 + A_2x_2 =b , \, x_1\in {\cal X}_1, x_2\in {\cal X}_2 \}.
   \]
Again, without ambiguity, some letters and notation are reused.

\subsection{Algorithm}

An alternative splitting version of the balanced ALM (\ref{SD-PPA1-P}) for the specific model (\ref{Problem-B2}) can be presented as below.

{\hdp{
{\bf Algorithm: An alternative splitting version of the balanced ALM for (\ref{Problem-B2})}
\smallskip

Let $r>0$, $s>0$ and $\delta>0$ be arbitrary constants. Define
\[\label{H2}
H_2=\frac{1}{s}A_2A_2^T + (\frac{1}{r} +\delta) I_m,
\]
$$
q_2^k:=x_2^k + \frac{1}{s}A_2^T\lambda^k \;\hbox{and}\; s_2^k=A_1(2{x}_1^{k+1}-x_1^k)+ A_2(2{x}_2^{k+1}-x_2^k)-b.
$$
Then, with $w^k= ({x}_1^{k},x_2^{k}, \lambda^{k})$, the new iterate ${w}^{k+1} =({x}_1^{k+1}, {x}_2^{k+1},{\lambda}^{k+1})$ is generated via the following steps:
\begin{subequations} \label{SD-PPA2-P}
\begin{numcases}{\hskip-1.2cm\;}
\label{SD2-xP}    {x}_1^{k+1}    = \arg\min \bigl\{\theta_1(x_1) - x_1^TA_1^T\lambda^k  +
                                 {\displaystyle\frac{r}{2}}
                                                                                 \|A_1(x_1-x_1^k)\|^2  +   {\displaystyle\frac{\delta}{2}}
                                                                                 \|x_1-x_1^k\|^2    \; \big| \; x_1\in{\cal X}_1 \bigr\},  \\
 \label{SD2-yP}    {x}_2^{k+1}    = \arg\min \bigl\{\theta_2(x_2)   +   {\displaystyle\frac{s}{2}}
                                                                                 \|x_2-q_2^k\|^2    \; \big| \; x_2\in{\cal X}_2 \bigr\},  \\
\label{SD2-lP}    {\lambda}^{k+1} =  \arg\min \Bigl\{ \frac{1}{2}(\lambda-\lambda^k)H_2(\lambda-\lambda_k)+(s_2^k)^T\lambda \;| \;\lambda\in \Lambda\Bigr\}.
\end{numcases}
\end{subequations}
}}

\begin{remark}
In the algorithm (\ref{SD-PPA2-P}), we see that only the $x_2$-subproblem (\ref{SD2-yP}) can be reduced to estimating the proximity operator of $\theta_2$ if ${\cal X}_2=\Re^{n_2}$, while the $x_1$-subproblem (\ref{SD2-xP}) is not proximity-induced because the term $\|A_1(x_1-x_1^k)\|^2$ is kept. As a balance, the matrix $H_2$ defined in (\ref{H2}) which determines the quadratic programming problem (\ref{SD2-lP}) does not involve $A_1$. In this sense, the $x_i$-subproblems and the $\lambda$-subproblem are balanced in another way. For the generic model (\ref{Problem-m}) with $p>2$, the splitter version of the balanced ALM (\ref{M-PRE}) can be revised in the sense that some of its $x_i$-subproblems are flexibly chosen to keep the terms $\|A_i(x_i-x_i^k)\|^2$ whilst the quadratic term determining the $\lambda^k$-subproblem does not involve the corresponding $A_i$'s. Thus, different algorithms with different balanced subproblems can be designed analogously. The algorithm (\ref{SD-PPA2-P}) is just the simplest illustration with $p=2$ for this philosophy.
\end{remark}

\subsection{Convergence results}

As mentioned, to prove the convergence of the algorithm (\ref{SD-PPA2-P}), we just need to prove an inequality similar as (\ref{PPA-T1F}) in Theorem \ref{HeB1} and show that the accompanying matrix is positive definite.

\begin{proposition}
Let $r>0$, $s>0$, and $\delta>0$ be arbitrary constants. Then, the matrix defined as
\[ \label{PPA-T2H}
  H  = \left(\begin{array}{ccc}
         r A_1^TA_1  +    \delta I_{n_1}    & 0  &      A_1^T \\
         0  &    s I_{n_2}      &      A_2^T \\
            A_1      & A_2     &    \displaystyle{\frac{1}{s}A_2A_2^T +  (\frac{1}{r} +\delta)} I_m
            \end{array}\right)\]
is positive definite.
\end{proposition}
\noindent{\bf Proof}.
Note that
$$    H=    \left(\begin{array}{ccc}
          r A_1^TA_1  +    \delta I_{n_1}   &0   &      A_1^T \\
          0    &0   &      0\\
            A      &0           &    \displaystyle{\frac{1}{r}}I_m
            \end{array}\right)   +
              \left(\begin{array}{ccc}
             0 & 0& 0\\
             &    s I_{n_2}      &     A_2^T \\
             0      &   A_2       &    \displaystyle{\frac{1}{s}}    A_2A_2^T
            \end{array}\right)   +  \left(\begin{array}{ccc}
           0  &   0    &     0 \\
           0  &   0    &     0 \\
          0   &   0    &   \delta I_m
            \end{array}\right).$$
For any $w=(x,y, \lambda)\ne 0$, we have
   $$    w^THw=  \Bigl(\Big\|\sqrt{r}A_1x+ \textstyle{\sqrt{\frac{1}{r}}}\lambda\Big\|^2 + \delta\|x\|^2 \Bigr)+   \Big\|\sqrt{s}y+ \textstyle{\sqrt{\frac{1}{s}}}A_2^T\lambda\Big\|^2  + \delta \|\lambda\|^2>0.  $$
    Thus,  the matrix $H$ is positive definite.  \hfill {$\Box$}

\begin{theorem} \label{HeB2} Let $\{{w}^{k}=(x_1^k,x_2^k,\lambda^k)\}$ be the sequence generated by the algorithm (\ref{SD-PPA2-P}) and $H$ be defined in (\ref{PPA-T2H}). Then, we have
\[  \label{PPA-T2F}  {w}^{k+1} \in \Omega, \;\;
      \theta(u)
      -\theta({u}^{k+1})  +(w-  {w}^{k+1})^TF({w}^{k+1}) \ge  (w- {w}^{k+1})^T H(w^k- {w}^{k+1}) , \;\;  \forall\, w\in
         \Omega.
          \]
      \end{theorem}
\noindent{\bf Proof}.  According to Lemma \ref{CP-TF},  ${x}^{k+1}$ generated by  \eqref{SD2-xP}  is characterized by the VI
$$ {x}_1^{k+1}\in {\cal X}_1,  \;\;
   \theta_1(x_1) - \theta_1({x}_1^{k+1}) +
 (x_1-x_1^{k+1})^T\bigl\{-A_1^T\lambda^k  + (r A_1^TA_1 +\delta) ({x}_1^{k+1} - x_1^k)  \bigr\}\ge 0,  \;\; \forall\, x_1\in {\cal X}_1.
 $$
Then, for any unknown ${\lambda}^{k+1}$, we have
\begin{eqnarray} \label{SD2-Pre-x}
 \lefteqn{ {x}_1^{k+1}\in {\cal X}_1,  \;\;
   \theta_1(x_1) - \theta_2({x}_1^{k+1}) +
 (x_1-{x}_1^{k+1})^T(-A_1^T{\lambda}^{k+1} )  } \nn \\
   & & \qquad \qquad \ge  (x_1-{x}_1^{k+1})^T \!\bigl\{ (r A_1^T\!A_1 +\delta) (x_1^k - {x}_1^{k+1})  + A_1^T\!(\lambda^k - {\lambda}^{k+1}) \bigr\},  \;\; \forall\, x_1\in {\cal X}_1.
   \end{eqnarray}
Analogously, it follows from Lemma \ref{CP-TF} that ${x}_2^{k+1}$ generated by  \eqref{SD2-yP}  can be characterized by the VI
$$ {x}_2^{k+1}\in {\cal X}_2,  \;\;
   \theta_2(x_2) - \theta_2({x}_2^{k+1}) +
 (x_2-{x}_2^{k+1})^T\bigl\{-A_2^T\lambda^k  + s ({x}_2^{k+1} - x_2^k)  \bigr\}\ge 0,  \;\; \forall\, x\in {\cal X}_2.
 $$
Then, for any unknown ${\lambda}^{k+1}$, we have
\begin{eqnarray} \label{SD2-Pre-y}
 \lefteqn{ {x}_2^{k+1}\in {\cal X}_2,  \;\;
   \theta_2(x_2) - \theta_2({x}_2^{k+1}) +
 (x_2-{x}_2^{k+1})^T(-A_2^T{\lambda}^{k+1} )  } \nn \\
   & & \qquad \qquad\qquad  \ge  (x_2-{x}_2^{k+1})^T  \bigl\{ s (x_2^k - {x}_2^{k+1})  + A_2^T(\lambda^k - {\lambda}^{k+1}) \bigr\},  \;\; \forall\, x_2\in {\cal X}_2.
   \end{eqnarray}
Similarly, according to Lemma \ref{CP-TF},  ${\lambda}^{k+1}$ generated by  \eqref{SD2-lP}  is characterized by the VI: Finding ${\lambda}^{k+1}\in \Lambda$ such that
      $$
(\lambda-{\lambda}^{k+1})^T\Bigl\{ \Bigl( A_1[2{x}_1^{k+1}- x_1^k] +A_2[2{x}_2^{k+1}- x_2^k]-b\Bigr)  +
   \Bigl(\frac{1}{s}A_2A_2^T +\bigl(\frac{1}{r} +\delta\bigr) I_m  \Bigr) ({\lambda}^{k+1}-\lambda^k) \Bigr\}
     \ge 0, \;\; \forall\, \lambda\in \Lambda.  $$
It can be rewritten as
\begin{eqnarray} \label{SD2-Pre-l}
 \hskip -0.8cm \lefteqn{
{\lambda}^{k+1}\in \Lambda, \quad  (\lambda-{\lambda}^{k+1})^T(A_1{x}_1^{k+1} +A_2{x}_2^{k+1} -b)  } \nn \\
  \hskip -0.3cm &\ge &\!\!\!  (\lambda-{\lambda}^{k+1})^T  \Bigl\{A_1(x_1^k - {x}_1^{k+1}) \! +\! A_2(x_2^k-{x}_2^{k+1})   \!+\!
      \Bigl(\displaystyle{\frac{1}{s}}A_2A_2^T \!\!
     +\! \bigl(\frac{1}{r} +\delta\bigr) I_m \Bigr) (\lambda^k - {\lambda}^{k+1}) \Bigr\},
   \end{eqnarray}
for all $\lambda\in \Lambda$. Combining  \eqref{SD2-Pre-x}, \eqref{SD2-Pre-y}  and  \eqref{SD2-Pre-l}, and  using the notation in \eqref{VI1}, we get the following assertion.
 \hfill {$\Box$}

\subsection{Comparison with linearized versions of the ADMM}

It is interesting to compare the proposed algorithm (\ref{SD-PPA2-P}) with the well-known linearized versions of the ADMM. For the model (\ref{Problem-B2}), the original ADMM scheme reads as
\begin{subequations} \label{ADMM-s}
\begin{numcases}{\quad}
\label{ADMMs-x}  x_1^{k+1} \in \arg\min \bigl\{ \theta_1(x_1) - x_1^TA^T\lambda^k +{\textstyle{ \frac{r}{2}}}\|A_1x_1+A_2x_2^k-b \|^2 \; \big| \; x_1\in{\cal X}_1
    \bigr\},\\[0.1cm]
\label{ADMMs-y}  x_2^{k+1}\in \arg\min \bigl\{\theta_2(x_2)  - {x}_2^TA_2^T\lambda^k
    +{\textstyle{ \frac{r}{2}}}\|A_1x_1^{k+1}+A_2x_2-b \|^2 \;   \; \big| \; x_2\in{\cal X}_2
    \bigr\},\\[0.1cm]
\label{ADMMs-l}{\lambda}^{k+1} = \lambda^k - r(A_1 x_1^{k+1} + A_2x_2^{k+1} -b),
\end{numcases}
\end{subequations}
in which $r>0$ is the penalty parameter and $\lambda\in \Re^m$ is the Lagrange multiplier. The first proximal version of the ADMM (PADMM) which suggests regularizing both the $x_1$- and $x_2$-subproblems in (\ref{ADMM-s}) with generic proximal terms was proposed in \cite{He02MP} (see also \cite{Eck94} for a special case). For simplicity, let us assume that the $x_1$-subproblem (\ref{ADMMs-x}) is easy but the $x_2$-subproblem (\ref{ADMMs-y}) is difficult. Then, the PADMM in \cite{He02MP} can be written as
{\small
\begin{subequations} \label{LADMM}
\begin{numcases}{\hskip-0.8cm \quad}
\label{LADMM-x}  x_1^{k+1} \in \arg\min \bigl\{ \theta_1(x_1) - x_1^TA_1^T\lambda^k +{\textstyle{ \frac{r}{2}}}\|A_1x_1+A_2x_2^k-b \|^2 \; \big| \; x_1\in{\cal X}_1
    \bigr\},\\[0.1cm]
\label{LADMM-y}  x_2^{k+1}\in \arg\min \bigl\{\theta_2(x_2)  - {x}_2^TA_2^T\lambda^k
    +{\textstyle{ \frac{r}{2}}}\|A_1x_1^{k+1}+A_2x_2-b \|^2 +{\textstyle{\frac{1}{2}}}\|x_2-x_2^k\|_G^2\;  \big| \; x_2\in{\cal X}_2
    \bigr\},\\[0.1cm]
\label{LADMM-l}{\lambda}^{k+1} = \lambda^k - r (A_1 x_1^{k+1} + A_2x_2^{k+1} -b),
\end{numcases}
\end{subequations}
}
\hskip-1mm in which $G \in \Re^{n_2\times n_2}$ is a positive definite matrix. Because of the same reason as mentioned for (\ref{PALM-E}), it is interesting to consider ``linearizing" the quadratic term ``$\frac{r}{2}\|A_1x_1^{k+1}+A_2x_2-b \|^2$" and thus alleviating the subproblem (\ref{LADMM-y}) as estimating the proximity operator of $\theta_2(x_2)$ when ${\cal X}_2=\Re^{n_2}$. Similar as (\ref{L1PALM-E}), this can be done by choosing $G:=sI_{n_2}-rA_2^TA_2$ in (\ref{LADMM-y}). As well discussed in the literature, e.g., \cite{Esser,Lin,COAP-HMY,WY,YY13}, for various applications arising in image processing, statistical learning, and others, the condition
   \[  \label{sYaoqiu}  s> r  \|A_2^TA_2\|
       \]
is required to ensure the positive definiteness of $G$ and thus the convergence of (\ref{LADMM}). Recently, the condition (\ref{sYaoqiu}) is further optimally improved in \cite{COAP-HMY} as $ s> 0.75 \cdot r \|A_2^TA_2\| $. Similarly as (\ref{PALM-E-x}) and (\ref{CP-x}), though $\theta_2(x_2)$ and $A_2$ are decoupled in notation if $G:=sI_{n_2}-rA_2^TA_2$ in (\ref{LADMM-y}), the subproblem (\ref{LADMM-y}) is correlated implicitly with $A_2$ via the condition (\ref{sYaoqiu}) or its improved one in \cite{COAP-HMY}. Hence, efficiency of all existing linearized versions of the ADMM is severely affected if $\|A_2^TA_2\|$ is large. In this sense, the algorithm (\ref{SD-PPA2-P}) improves existing linearized versions of the ADMM in the sense that the $x_2$-subproblem (\ref{SD2-yP}) can be reduced to estimating the proximity operator of $\theta_2$ if ${\cal X}_2=\Re^{n_2}$, while it is not affected by $\|A_2^TA_2\|$ and thus possible tiny step sizes could be avoided even if $\|A_2^TA_2\|$ is large.

\section{More generalized versions}\label{Sec:remarks}

\setcounter{equation}{0}

In the preceding sections, the balanced ALM (\ref{SD-PPA1-P}) is proposed for the generic model (\ref{Problem-B1}), and then its splitting versions (\ref{M-PRE}) and (\ref{SD-PPA2-P}) are studied for the separable models (\ref{Problem-m}) and (\ref{Problem-B2}), respectively. As mentioned, it was shown in \cite{Rock76A} that the classic ALM (\ref{ALM-E}) is an application of the PPA proposed in \cite{Mar70}. In view of the generalized version of the PPA studied in \cite{GT}, all the proposed algorithms (\ref{SD-PPA1-P}), (\ref{M-PRE}) and (\ref{SD-PPA2-P}) can be further generalized. For instance, the balanced ALM (\ref{SD-PPA1-P}) can be generalized as
\begin{subequations} \label{GSD-PPA1-P}
\begin{numcases}{\hskip-1.5cm\;}
\label{GSD-xP}   {\tilde x}^{k} =
  \arg\min  \bigl\{   \theta(x)  + \frac{r}{2}\|x- q_0^k\|^2  \; |\;  x\in {\cal X} \bigr\},  \\
\label{GSD-yP}    {\tilde \lambda}^{k} =  \arg\min \bigl\{
    \frac{1}{2}(\lambda-\lambda^k)H_0(\lambda-\lambda^k)+({\tilde{s}}_0^k)^T{\lambda} \;|\; \lambda\in \Lambda\bigr\},\\
\label{GSD-relax} \left(\begin{array}{c}
                    x^{k+1}\\
                    \lambda^{k+1}
                  \end{array}\right)=\left(\begin{array}{c}
                    x^{k}\\
                    \lambda^{k}
                  \end{array}\right)
                  -\alpha\left(\begin{array}{c}
                    x^{k}-{\tilde x}^k\\
                    \lambda^{k}-{\tilde \lambda}^k
                  \end{array}\right) \;\hbox{with} \;\alpha \in (0,2),
\end{numcases}
\end{subequations}
where $\tilde{s}_0^k = A(2\tilde{x}^k - x^k) -b$. Clearly, the scheme (\ref{GSD-PPA1-P}) includes the balanced ALM (\ref{SD-PPA1-P}) as a special case with $\alpha=1$. Numerically, the extra step (\ref{GSD-relax}) has been shown to be able to accelerate the convergence of the classic PPA for various problems. We refer to, e.g., \cite{Ber82,COAP-HYZ,HY-SIIMS}, for some empirical studies. Hence, it is motivated to consider the generalized scheme (\ref{GSD-PPA1-P}) to replace the balanced ALM (\ref{SD-PPA1-P}).

To establish the convergence of (\ref{GSD-PPA1-P}), we just need to follow the roadmap in Section \ref{Sec:convergence} and prove some similar theorems. For instance, the inequality (\ref{Hau1-A0}) in Theorem \ref{HB1-A} can be generalized as
 \begin{eqnarray*} 
 \lefteqn{  \alpha\bigl( {\theta}(x) - {\theta}(\tilde{x}^k) +(w- \tilde{w}^k)^T
             F(w) \bigr) }    \nn \\
             &\ge & \frac{1}{2}
    \bigl(\|w-w^{k+1}\|_{H}^2 -\|w-w^k\|_{H}^2  \bigr) +
        \frac{\alpha(2-\alpha)}{2}\|w^k-\tilde{w}^k\|_{H}^2, \;\;  \forall w\in \Omega.
        \end{eqnarray*}
Moreover, the inequality (\ref{Hau1-Bw}) in Theorem \ref{HB1-B} can be generalized as
$$
 \|w^{k+1}-w^*\|_H^2  \le \|w^k-w^* \|_H^2  -\alpha(2-\alpha)\|w^k-\tilde{w}^k\|_{H}^2, \;\;
 \forall w^*\in \Omega^*.
$$
Then, based on these inequalities, analogous as the analysis in Section \ref{Sec:convergence}, convergence results for the generalized version of the balanced ALM (\ref{GSD-PPA1-P}) can be obtained trivially.

In addition, the extra step (\ref{GSD-relax}) can be combined with the splitting versions of the balanced ALM (\ref{M-PRE}) and (\ref{SD-PPA2-P}) as well, and thus some generalized versions of the algorithms (\ref{M-PRE}) and (\ref{SD-PPA2-P}) can also be proposed. The details are omitted for succinctness.

\section{Conclusions}\label{Sec:conclusion}

In this paper, we reshape the classic augmented Lagrangian method (ALM) by balancing its subproblems. Convex programming problems with both linear equality and inequality constraints are considered. We propose a balanced ALM for the generic case, and various splitting versions for the separable cases. The balanced ALM and its splitting versions have the common feature that the subproblems are better balanced, and they are easier to be implemented for various applications. The balanced ALM advances the classic ALM by enlarging its applicable range, better balancing its subproblems, and improving its implementation. The balanced ALM and its splitting versions substantially enhance the rich literature of the classic ALM and its variants from a novel perspective, and open up the door to designing other application-tailored algorithms of the same kind for more specific/complicated problems.


{
 }

\end{CJK*}
\end{document}